\documentclass[draft,11pt,english]{article}
\usepackage{amsfonts,amsmath,amsthm,amscd,amssymb,latexsym}

\usepackage[T2A]{fontenc}
\usepackage[cp1251]{inputenc}
\usepackage[russian, ukrainian]{babel}

\def\XXint#1#2#3{{\setbox0=\hbox{$#1{#2#3}{\int}$}
\vcenter{\hbox{$#2#3$}}\kern-.5\wd0}}

\begin{document}

\title{On the global behavior of solutions of the Beltrami equations}
\author{Ruslan  Salimov, Mariia  Stefanchuk}
\date{}
\theoremstyle{plain}
\newtheorem{theorem}{Теорема}[section]
\newtheorem{lemma}{Лема}[section]
\newtheorem{proposition}{Твердження}[section]
\newtheorem{corollary}{Наслідок}[section]
\newtheorem{definition}{Означення}[section]
\theoremstyle{definition}
\newtheorem{example}{Приклад}[section]
\newtheorem{remark}{Зауваження}[section]
\newcommand{\keywords}{\textbf{Keywords. }\medskip}
\newcommand{\subjclass}{\textbf{MSC 2020. }\medskip}
\renewcommand{\abstract}{\textbf{Abstract.}\medskip}
\renewcommand{\refname}{References}
\numberwithin{equation}{section}

\maketitle

\begin{abstract}
In this paper, the estimate for growth of homeomorphic solutions of the Beltrami equation  at infinity is obtained, provided that the dilatation quotient has a global finite mean oscillation.

\end{abstract}
\medskip



\vskip10pt

\section{Introduction}

Let $D$ be a domain in the complex plane ${\Bbb C}$, i.e., a
connected and open subset of ${\Bbb C}$, and let $\mu:D\to{\Bbb C}$
be a measurable function with $|\mu(z)|<1$ a.e. (almost everywhere)
in $D$. The {\it Beltrami equation} is the equation of the form
\begin{equation}\label{eq1.1} f_{\overline{z}}\,=\,\mu(z)
f_z\,\end{equation} where
$f_{\overline{z}}={\overline{\partial}}f=(f_x+if_y)/2$,
$f_{z}=\partial f=(f_x-if_y)/2$, $z=x+iy$, and $f_x$ and $f_y$ are
partial derivatives of $f$ in $x$ and $y$, correspondingly. The
function $\mu$ is called the {\it complex coefficient} and
\begin{equation}\label{eq1.1a}K_{\mu}(z)\,=\,\frac{1+|\mu(z)|}{1-|\mu (z)|}\end{equation}
the {\it dilatation quotient} for the equation (\ref{eq1.1}). The
Beltrami equation (\ref{eq1.1}) is said to be {\it degenerate} if
${\rm ess}\,{\rm sup}\,K_{\mu}(z)=\infty$. The existence theorem for
homeomorphic $W^{1,1}_{\mathrm{loc}}$ solutions was established to
many degenerate Beltrami equations, see, e.g., related references in
the recent monographs  \cite{AIM}, \cite{MRSY},  \cite{GRSY1}; cf. also \cite{GRSY2}, \cite{RSY$_0$}~--~\cite{RSY$_4$}.

Recall that the {\it (conformal) modulus} of a family $\Gamma$ of
curves $\gamma$ in ${\Bbb C}$ is the quantity
\begin{equation}\label{eq4132}
M(\Gamma)=\inf_{\rho \in \,{\rm adm}\,\Gamma} \int\limits_{{\Bbb C}}
\rho ^2 (z)\ \ dx\,dy\,
\end{equation}
where a Borel function $\rho:{\Bbb C}\,\rightarrow [0,\infty]$ is {\it
admissible} for $\Gamma$, write $\rho \in {\rm adm} \,\Gamma $, if
\begin{equation}\label{eq4133}
\int\limits_{\gamma}\rho \,\,ds\geqslant 1\ \ \ \ \ \ \ \forall\ \gamma
\in \Gamma\,
\end{equation}
where $s$ is a natural parameter of the length on $\gamma$.

Throughout this paper, $$B(z_{0},\,r)=\{ z\in{\Bbb C}: |z-z_{0}|<r\}\,,$$
 $$S(z_{0},\,r)=\{ z\in{\Bbb C}: |z-z_{0}|=r\}\,, $$
 and
 $$\mathbb{A}(z_0, r_1,r_2)=\{ z\,\in\,{\Bbb C} : r_1<|z-z_0|<r_2\}\,.$$
 Let $E,$ $F\subset\overline{{\Bbb C}}$ be arbitrary
sets. Denote by $\Delta(E,F,D)$  a family of all curves
$\gamma:[a,b]\rightarrow \overline{{\Bbb C}}$ joining $E$ and $F$ in
$D,$ i.e., $\gamma(a)\in E,\gamma(b) \in F$ and $\gamma(t)\in D$ as
$t \in (a, b).$

Here a {\it condenser} is a pair $\mathcal{E}=(A,C)$ where
$A\subset \Bbb C$ is open and $C$ is a non-empty compact set
contained in $A$. $\mathcal{E}$ is a {\it  ringlike condenser}  if $B=A\setminus C$ is a ring, i.e., if $B$ is a domain
whose complement $\overline{\mathbb{C}}\setminus B$ has exactly two components where
$\overline{\mathbb{C}}=\mathbb{C}\cup\{\infty\}$ is the one-point compactification of $\mathbb{C}$.
$\mathcal{E}$ is a {\it  bounded condenser}  if $A$ is
bounded. A condenser $\mathcal{E}=(A, C)$ is said to be in a domain $G$ if
$A\subset G$.

\medskip

The following lemma is immediate.

\vskip 2mm

\textbf{Lemma 1.} \label{lem33.1} {\it If $f:G\to{\Bbb C}$ is open and $\mathcal{E}=(A,C)$ is a
condenser in $G$, then $\left(fA,fC\right)$ is a condenser in $fG$.}

\vskip 2mm

In the above situation we denote $f\mathcal{E}=\left(fA,fC\right)$.

\medskip

Let $\mathcal{E}=\left(A,C\right)$ be a condenser. We set
$$\mathrm{cap}\,\mathcal{E}=\mathrm{cap}\left(A,C\right)=\inf
\limits_{u\in W_0\left(\mathcal{E}\right)}\int\limits_{A}\vert\nabla u\vert^2\
dxdy $$ and call it the {\it  capacity} of the condenser
$\mathcal{E}\,.$ The set $W_0(\mathcal{E})=W_0(A,C)$ is the family of nonnegative
functions $u:A\to \mathbb{R}$ such that $u\in C_0(A)$,  $u(z)\geqslant 1$
for $z\in C$, and  $u$ is absolutely continuous on  lines ($\textrm{ACL}$). In the above formula
$$\vert\nabla u\vert= \sqrt{\left(\frac{\partial u}{\partial x}\right)^2+\left(\frac{\partial u}{\partial y}\right)^2.}$$

We mention some properties of the capacity of a condenser. It was proven in (\cite{Sh}, Theorem 1)  that
\begin{equation}\label{Shlyk}
\mathrm{cap}\,\mathcal{E}=M(\Delta(\partial A,\partial C;A\setminus C ))\,,\end{equation}
where $\Delta(\partial A,\partial C;A\setminus C)$  denotes the set of all continuous curves joining the
boundaries $\partial A$ and $\partial C$ in $A\setminus C$.

\vskip 2mm

Moreover, the following estimate is known: \begin{equation}\label{mazn} {\rm cap}\,\mathcal{E}\geqslant \frac{4 \pi }{\log \,  \frac{m(A)}{m(C)}}\end{equation}
(see, e.g., (8.8) in  \cite{Maz}).

\vskip 2mm

The following notion is motivated by the ring definition of Gehring
for qua\-si\-con\-for\-mal mappings, see, e.g., \cite{Ge$_3$},
introduced first in the plane, see \cite{RSY$_3$}, and extended
later on to the space case in \cite{RS}, see also Chapters 7 and 11
in \cite{MRSY}, cf. \cite{AC$_1$}, \cite{ACS}, \cite{Cr3}, \cite{RS_{2}}.

Given a domain $D$ in $\mathbb{C}$, a (Lebesgue) measurable function
$Q:D\rightarrow\,[0,\infty]$, $z_0\in D,$
a homeomorphism $f:D\rightarrow \overline{\mathbb{C}}$ is said to be a
{\it ring $Q$--homeomorphism at the point $z_0$} if
\begin{equation}\label{eq1}
M\left(f\left(\Delta\left(S_1,\,S_2,\,\mathbb{A}(z_0, r_1,r_2)\right)\right)\right)\
\leqslant \int\limits_{\mathbb{A}(z_0, r_1,r_2)} Q(z)\cdot \eta^2(|z-z_0|)\ dx\,dy
\end{equation}
for every ring $\mathbb{A}(z_0, r_1,r_2)$ and the circles $S_i=S(z_0, r_i)$, $i=1,2,$
where $0<r_1<r_2< r_0\,\colon =\,{\rm dist}\, (z_0,\partial D),$ and
every measurable function $\eta : (r_1,r_2)\rightarrow [0,\infty ]$
such that
$$\int\limits_{r_1}^{r_2}\eta(r)\ dr \ =\ 1\,.$$
$f$ is called a {\it ring $Q$--homeomorphism in the domain} $D$ if
$f$ is a ring $Q$--ho\-meo\-mor\-phism at every point $z_0\in D$.

The following statement was first proved in \cite{LSS}, Theorem 3.1,  cf. also Corollory  3.1 in \cite {Sa}.

\vskip 2mm

\textbf{Proposition 1.} \label{th3.1989898} {\it Let $f$ be a homeomorphic $W^{1,1}_{\rm loc}$ solution of the
Beltra\-mi equation (\ref{eq1.1}). Then $f$ is a ring
$Q$-homeomorphism at each point $z_0\in D$ with $Q(z)=K_{\mu}(z)$.}

\vskip 2mm

\section{GFMO functions}

Similarly to \cite{IR} (cf. also  \cite{RSY$_0$},  \cite{RSY$_2$}), we say that a function $\varphi\colon\mathbb{C}\to\mathbb{R}$ has \emph{global finite mean oscillation at a point} $z_{0}\in\mathbb{C}$, abbr. $\varphi\in GFMO(z_{0})$, if
\begin{equation}\label{11}
\limsup\limits_{R\to\infty}\frac{1}{m\left(B(z_{0},R)\right)}\int\limits_{B(z_{0},R)}|\varphi(z)-\overline{\varphi}_{R}|\,dxdy<\infty,
\end{equation}
where
$$\overline{\varphi}_{R}=\frac{1}{m\left(B(z_{0},R)\right)}\int\limits_{B(z_{0},R)}\varphi(z)\,dxdy$$
is the mean value of the function $\varphi(z)$ over $B(z_{0},R)$, $R>0$. Here $B(z_{0},R)=\{z\in\mathbb{C}\colon|z-z_{0}|<R\}$, and condition \eqref{11} includes the assumption that $\varphi$ is integrable in $B(z_{0},R)$ for $R>0$.

\vskip 2mm

\textbf{Proposition 2.} {\it If, for some collection of numbers $\varphi_{R}\in\mathbb{R}$, $R\in[r_{0},+\infty)$, $r_{0}>0$,
$$\limsup\limits_{R\to\infty}\frac{1}{m\left(B(z_{0},R)\right)}\int\limits_{B(z_{0},R)}|\varphi(z)-\varphi_{R}|\,dxdy<\infty,$$
then $\varphi$ has global finite mean oscillation at $z_{0}$.}

\vskip 2mm

\textbf{Proof.} Indeed, by the triangle inequality,
$$\frac{1}{m\left(B(z_{0},R)\right)}\int\limits_{B(z_{0},R)}|\varphi(z)-\overline{\varphi}_{R}|\,dxdy\leqslant$$
$$\leqslant\frac{1}{m\left(B(z_{0},R)\right)}\int\limits_{B(z_{0},R)}|\varphi(z)-\varphi_{R}|\,dxdy+|\varphi_{R}-\overline{\varphi}_{R}|\leqslant$$
$$\leqslant\frac{2}{m\left(B(z_{0},R)\right)}\int\limits_{B(z_{0},R)}|\varphi(z)-\varphi_{R}|\,dxdy.$$

\vskip 2mm

\textbf{Corollary 1.} {\it If, for a point $z_{0}\in\mathbb{C}$,
$$\limsup\limits_{R\to\infty}\frac{1}{m\left(B(z_{0},R)\right)}\int\limits_{B(z_{0},R)}|\varphi(z)-\varphi(z_{0})|\,dxdy<\infty,$$
then $\varphi$ has global finite mean oscillation at $z_{0}$.}

\vskip 2mm

\textbf{Corollary 2.} {\it If, for a point $z_{0}\in\mathbb{C}$,
$$\limsup\limits_{R\to\infty}\frac{1}{m\left(B(z_{0},R)\right)}\int\limits_{B(z_{0},R)}|\varphi(z)|\,dxdy<\infty,$$
then $\varphi$ has global finite mean oscillation at $z_{0}$.}

\vskip 2mm

\textbf{Lemma 2.} {\it Let $z_{0}\in\mathbb{C}$. If a nonnegative function $\varphi\colon\mathbb{C}\to\mathbb{R}$ has  global finite mean oscillation at $z_{0}$ and $\varphi$ is integrable in $B(z_{0},e)$, then, for $R>e^{e},$
$$\int\limits_{\mathbb{A}(z_{0},e,R)}\frac{\varphi(z)\,dxdy}{\left(|z-z_{0}|\log|z-z_{0}|\right)^{2}}\leqslant C\cdot\log\log R,$$
where
$$C=\frac{\pi}{6}((24+\pi^{2})e^2\delta_{\infty}+2\pi^2\varphi_{0}),$$
$\varphi_{0}$ is the mean value of $\varphi$ over the disk $B(z_{0},e)$ and $$\delta_{\infty}=\delta_{\infty}(\varphi)=\sup\limits_{R\in(e^{e},+\infty)}\frac{1}{m\left(B(z_{0},R)\right)}\int\limits_{B(z_{0},R)}|\varphi(z)-\overline{\varphi}_{R}|\,dxdy$$ is the maximal dispersion of $\varphi$.}

\vskip 2mm

\textbf{Proof.} Let $R>e^{e}$, $r_{k}=e^{k}$, $\mathbb{A}_{k}=\{z\in\mathbb{C}\colon r_{k}\leqslant|z-z_{0}|<r_{k+1}\}$.
Clearly,
$$\delta_{\infty}=\sup\limits_{R\in(e^{e},+\infty)}\frac{1}{m\left(B(z_{0},R)\right)}\int\limits_{B(z_{0},R)}|\varphi(z)-\overline{\varphi}_{R}|\,dxdy<\infty\,, $$
$B_{k}=B(z_{0},r_{k})$ and let $\varphi_{k}$ be the mean value of $\varphi(z)$ over $B_{k}$, $k=1,2,...$. Take a natural number $N$ such that $R\in[r_{N},r_{N+1})$.

Then $\mathbb{A}(z_{0},e,R)\subset\Delta(R)=\bigcup\limits_{k=1}^{N}\mathbb{A}_{k}$ and
$$I(R)=\int\limits_{\Delta(R)}\varphi(z)\alpha(|z-z_{0}|)\,dxdy\leqslant|S_{1}(R)|+S_{2}(R)\,,$$
$$\alpha(t)=\frac{1}{(t\log t)^{2}}\,,$$
$$S_{1}(R)=\sum\limits_{k=1}^{N}\int\limits_{\mathbb{A}_{k}}(\varphi(z)-\varphi_{k+1})\alpha(|z-z_{0}|)\,dxdy\,,$$
and
$$S_{2}(R)=\sum\limits_{k=1}^{N}\varphi_{k+1}\int\limits_{\mathbb{A}_{k}}\alpha(|z-z_{0}|)\,dxdy\,.$$
Since $\mathbb{A}_{k}\subset B_{k+1}$, $\frac{1}{|z-z_{0}|^{2}}\leqslant\frac{\pi e^{2}}{m(B_{k+1})}$ for $z\in\mathbb{A}_{k}$
and $\log|z-z_{0}|>k$ in $\mathbb{A}_{k}$, then
$$|S_{1}(R)|\leqslant\pi e^{2}\sum\limits_{k=1}^{N}\frac{1}{k^{2}}\cdot\frac{1}{m(B_{k+1})}\int\limits_{B_{k+1}}|\varphi(z)-\varphi_{k+1}|\,dxdy\leqslant$$
$$\leqslant\pi e^{2}\delta_{\infty}\sum\limits_{k=1}^{N}\frac{1}{k^{2}}\leqslant \pi e^{2}\delta_{\infty}\sum\limits_{k=1}^{\infty}\frac{1}{k^{2}}=
\frac{\pi^{3}e^{2}\delta_{\infty}}{6}\,.$$
Now,
$$\int\limits_{\mathbb{A}_{k}}\alpha(|z-z_{0}|)\,dxdy\leqslant\frac{1}{k^{2}}\int\limits_{\mathbb{A}_{k}}\frac{dxdy}{|z-z_{0}|^{2}}=\frac{2\pi}{k^{2}}\,.$$
Moreover,
$$|\varphi_{k-1}-\varphi_{k}|=\left|\frac{1}{m(B_{k-1})}\int\limits_{B_{k-1}}\varphi(z)\,dxdy-\frac{1}{m(B_{k-1})}\int\limits_{B_{k-1}}\varphi_{k}\,dxdy\right|\leqslant$$
$$\leqslant\frac{1}{m(B_{k-1})}\int\limits_{B_{k-1}}|\varphi(z)-\varphi_{k}|\,dxdy\leqslant
\frac{e^{2}}{m(B_{k})}\int\limits_{B_{k}}|\varphi(z)-\varphi_{k}|\,dxdy\leqslant e^{2}\delta_{\infty}\,,$$
and by the triangle inequality, for $k\geqslant1$
$$\varphi_{k+1}=|\varphi_{k+1}|=\left|\varphi_{1}+\sum\limits_{l=2}^{k+1}(\varphi_{l}-\varphi_{l-1})\right|\leqslant$$
$$\leqslant|\varphi_{1}|+\sum\limits_{l=2}^{k+1}|\varphi_{l}-\varphi_{l-1}|\leqslant |\varphi_{1}|+e^{2}\delta_{\infty}\, k.$$
Hence,
$$S_{2}(R)=|S_{2}(R)|\leqslant2\pi\sum\limits_{k=1}^{N}\frac{\varphi_{k+1}}{k^{2}}\leqslant
2\pi\sum\limits_{k=1}^{N}\frac{\varphi_{1}+e^{2}\delta_{\infty}\, k}{k^{2}}\leqslant$$
$$\leqslant2\pi\varphi_{1}\sum\limits_{k=1}^{\infty}\frac{1}{k^{2}}+2\pi e^{2}\delta_{\infty} \, \sum\limits_{k=1}^{N}\frac{1}{k}=$$
$$=\frac{\pi^{3}\varphi_{1}}{3}+2\pi e^{2}\delta_{\infty}\, \sum\limits_{k=1}^{N}\frac{1}{k}\,.$$
But
$$\sum\limits_{k=2}^{N}\frac{1}{k}<\int\limits_{1}^{N}\frac{dt}{t}=\log N$$
and, for $R>r_{N}$,
$$N=\log r_{N}<\log R\,.$$
Consequently,
$$\sum\limits_{k=1}^{N}\frac{1}{k}<1+\log\log R$$
and thus, for $R\in(e^{e},+\infty)$
$$I(R)\leqslant \frac{\pi^{3}e^{2}\delta_{\infty}}{6}+ \frac{\pi^{3}\varphi_{1}}{3}+2\pi e^{2}\delta_{0}(1+\log\log R)=$$
$$=\left(\frac{\pi^3e^2\delta_{\infty}+12\pi e^2\delta_{\infty} +2\pi^3\varphi_{1}}{6\log\log R}+2\pi e^{2}\delta_{\infty}\right)\log\log R\leqslant$$
$$\leqslant\frac{\pi}{6}((24+\pi^{2})e^2\delta_{\infty}+2\pi^2\varphi_{1})\log\log R\,.$$

Finally,
$$\int\limits_{\mathbb{A}(z_{0},e,R)}\frac{\varphi(z)\,dxdy}{\left(|z-z_{0}|\log|z-z_{0}|\right)^{2}}\leqslant I(R)\leqslant$$ $$\leqslant\frac{\pi}{6}((24+\pi^{2})e^2\delta_{\infty}+2\pi^2\varphi_{1})\log\log R\,.$$

\section{The behavior of homeomorphic so\-lutions of the Beltrami equations at infinity}

Set $$l_{f}(z_{0},e)=\min\limits_{|z-z_{0}|=e}|f(z)-f(z_{0})|\,,$$
$$\delta_{\infty}=\delta_{\infty}\left(K_{\mu},z_{0}\right)=$$
$$=\sup\limits_{R\in(e^{e},+\infty)}\frac{1}{m\left(B(z_{0},R)\right)}\int\limits_{B(z_{0},R)}|K_{\mu}(z)-K_{\mu,z_{0}}(R)|\,dxdy\,,$$ $$K_{\mu,z_{0}}(R)=\frac{1}{m\left(B(z_{0},R)\right)}\int\limits_{B(z_{0},R)}K_{\mu}(z)\,dxdy\,,\ \  k_0= K_{\mu,z_{0}}(e)\,.$$

\medskip

\textbf{Theorem 1.} {\it Let $\mu\colon \mathbb{C}\to \mathbb{C}$ be a measurable function with $|\mu(z)|<1$ a.e. and $f\colon \mathbb{C}\to \mathbb{C}$ be a homeomorphic $W_{{\rm loc}}^{1,1}$ solution of the Beltrami equation \eqref{eq1.1}. If $K_{\mu}\in GFMO(z_{0})$, $z_{0}\in\mathbb{C}$, then
\begin{equation}\label{664}
\liminf\limits_{R\to\infty}\frac{\max\limits_{|z-z_{0}|=R}|f(z)-f(z_{0})|}{(\log R)^{\frac{2\pi}{C}}}\geqslant l_{f}(z_{0},e)\,,
\end{equation}
where  $C=\frac{\pi}{6}((24+\pi^{2})e^{2}\delta_{\infty}+2\pi^{2}k_{0})$.}

\medskip

\textbf{Proof.} Consider the ring $\mathbb{A}(R)=\mathbb{A}(z_{0},e,R)$, with $R>e^{e}$. Set $\mathcal{E}=(B(z_{0},R),\,\overline{B(z_{0},e)})$. Then, by lemma 1 $f\mathcal{E}=(fB(z_{0},R),\,f\overline{B(z_{0},e)})$ is a condenser in $\mathbb{C}$, according \eqref{Shlyk},
$$\mathrm{cap}\,(fB(z_{0},R),\,f\overline{B(z_{0},e)})=M(\Delta(\partial f B(z_{0},e),\partial f B(z_{0},R);f\mathbb{A}(R)))$$
and, in view of the homeomorphism of $f$,
$$\Delta(\partial fB(z_{0},e),\partial fB(z_{0},R);f\mathbb{A}(R))=f\Delta(\partial B(z_{0},e),\partial B(z_{0},R);\mathbb{A}(R))\,.$$
By Proposition 1 $f$ is a ring $Q$-homeomorphism with $Q=K_{\mu}(z)$
\begin{equation}\label{66}
\mathrm{cap}\,(fB(z_{0},R),\,f\overline{B(z_{0},e)})\leqslant\int\limits_{\mathbb{A}(R)}K_{\mu}(z)\eta^{2}(|z-z_{0}|)\,dxdy
\end{equation}
for every measurable function $\eta \colon (e,R)\to [0,+\infty]$ such that $$\int\limits_{e}^{R}\eta(t)\,dt=1.$$ Choosing in \eqref{66} $\eta(t)=\frac{1}{t\log t\cdot\log\log R}$, we obtain
$$\mathrm{cap}\,(fB(z_{0},R),\,f\overline{B(z_{0},e)})\leqslant\frac{1}{(\log\log R)^{2}}\cdot\int\limits_{\mathbb{A}(R)}\frac{K_{\mu}(z)\,dxdy}{(|z-z_{0}|\log|z-z_{0}|)^{2}}\,.$$
Since $K_{\mu}\in GFMO(z_{0})$, then by lemma 2
\begin{equation}\label{661}
\mathrm{cap}\,(fB(z_{0},R),\,f\overline{B(z_{0},e)})\leqslant\frac{C}{\log\log R}\,,
\end{equation}
where $C=\frac{\pi}{6}((24+\pi^{2})e^{2}\delta_{\infty}+2\pi^{2}k_0).$
On the other hand, by \eqref{mazn}, we have
\begin{equation}\label{662}
\mathrm{cap}\,(fB(z_{0},R),\,f\overline{B(z_{0},e)})\geqslant\frac{4\pi}{\log\frac{m(fB(z_{0},R))}{m(f\overline{B(z_{0},e)})}}\,.
\end{equation}
Combining \eqref{661} and \eqref{662}, we obtain
$$\frac{4\pi}{\log\frac{m(fB(z_{0},R))}{m(f\overline{B(z_{0},e)})}}\leqslant\frac{C}{\log\log R}\,.$$
This gives
$$m(f\overline{B(z_{0},e)})\leqslant\frac{m(fB(z_{0},R))}{(\log R)^{\frac{4\pi}{C}}}\,.$$
Using the inequalities
$$\pi\left(\min\limits_{|z-z_{0}|=e}|f(z)-f(z_{0})|\right)^{2}\leqslant m(f\overline{B(z_{0},e)})\leqslant$$
$$\leqslant m(fB(z_{0},R))\leqslant
\pi\left(\max\limits_{|z-z_{0}|=R}|f(z)-f(z_{0})|\right)^{2}\,,$$
we obtain
\begin{equation}\label{663}
\min\limits_{|z-z_{0}|=e}|f(z)-f(z_{0})|\leqslant\frac{\max\limits_{|z-z_{0}|=R}|f(z)-f(z_{0})|}{(\log R)^{\frac{2\pi}{C}}}\,.
\end{equation}
Set
$$l_{f}(z_{0},e)=\min\limits_{|z-z_{0}|=e}|f(z)-f(z_{0})|\,.$$
Passing to the lower limit as $R\to\infty$ in \eqref{663}, we obtain relation \eqref{664}.

\newpage
Ruslan  Salimov,

Institute of Mathematics of the National Academy of Sciences of

Ukraine,

01601, Kyiv, Tereshchenkivska st. 3;

ruslan.salimov1@gmail.com

\medskip
Mariia  Stefanchuk,

Institute of Mathematics of the National Academy of Sciences of

Ukraine,

01601, Kyiv, Tereshchenkivska st. 3;

stefanmv43@gmail.com

\end{document}